\documentclass[10pt,reqno,a4paper]{amsart}
\usepackage{amsmath,amsfonts,amssymb,amsthm,amsxtra}

\newtheorem{thm}{Theorem}[section]

\newtheorem{prop}{Proposition}[section]
\newtheorem{exam}{Example}[section]
\newtheorem{defin}{Definition}[section]

\newcommand{\lc}{{[ \![}}
\newcommand{\rc}{{] \!]}}

\newcommand{\lL}{\mathbb{L}}

\newcommand{\T}{{\mathbb T}}

\newcommand{\R}{{\mathbb R}}

\begin{document}
\title{Poisson Algebras of Admissible Functions Associated to Twisted Dirac Structures}
\author{Alexander Cardona}
\address{Mathematics Department
\\Universidad de Los Andes
\\A.A. 4976 Bogot\'a, Colombia.}
\email{acardona@uniandes.edu.co}
\date{February 27, 2011}

\maketitle

\begin{abstract}
We define algebras of admissible functions associated to twisted
Dirac structures, and we show that they are Poisson algebras. We
study the standard cases associated to Dirac structures defined by
graphs of non-degenerate $2$-forms.
\end{abstract}

\vspace{0.5cm}



\vspace{0.5cm}

\section{Introduction}

Poisson algebras of admissible functions associated to non-twisted
Dirac structures have been studied by Courant and Weinstein (see
\cite{C}\cite{CW}). A Dirac structure on a manifold $M$ is a
maximally isotropic sub-bundle $\lL$ of the {\em generalized
tangent bundle} $\T M = TM \oplus T^* M$, which is involutive
under the {\em Courant bracket} on $\T M$
\begin{equation}
  \label{E:CB}
  [ X \oplus \xi, Y \oplus \eta ]_C = [X, Y] + \mathcal{L}_X \eta  -
\mathcal{L}_Y \xi - {1 \over 2}d (\iota_X \eta - \iota_Y \xi),
\end{equation}
where $X \oplus \xi, Y \oplus \eta \in \Gamma (\T M)$. The
isotropy condition here is given with respect to the natural
symmetric pairing
\begin{equation}\label{E:SPairing}
    \langle X \oplus \xi, Y \oplus \eta \rangle_{+} = {1 \over 2} (\xi (Y) + \eta (X))
\end{equation}
in $\T M$, and the bracket (\ref{E:CB}) can be seen as the
skew-symmetrization of the Dorfman bracket \cite{D}
\begin{equation}
 \label{E:DB}
  [ X \oplus \xi, Y \oplus \eta ]_D = [X, Y] + \mathcal{L}_X \eta  -
i_Y d \xi ,
\end{equation}
which coincides with (\ref{E:CB}) on sections of $\lL$. Particular
cases of Dirac manifolds are Poisson and symplectic manifolds
(which correspond to graphs, in the generalized tangent bundle $\T
M$, of the corresponding Poisson bi-vector and symplectic form,
respectively). In the symplectic case, for example, if $h \in
\Omega^2(M)$ denotes the symplectic form,
\begin{equation}\label{E:SymplGraph}
    \lL_h =  \{(X, i_X h) \in \Gamma (\T M)\mid X \in \frak{X}(M) \}
\end{equation}
defines a Dirac structure on $M$, and many features of the
symplectic geometry associated to $h$ are captured by this Dirac
structure. In particular, the Poisson algebra on $C^\infty (M)$
defined by the action of Hamiltonian vector fields on smooth
functions appears here as the algebra of admissible functions
associated to the Dirac structure (see e.g.\cite {C}\cite{CW}).
\\ \\
In general, given a Dirac structure $\lL$ on $M$, it is possible
to associate to it a Poisson algebra of smooth functions on $M$,
which is usually a subalgebra of $C^\infty (M)$, and is called the
algebra of \emph{admissible functions} \cite{C}. A smooth function
$f$ on a manifold $M$ with a Dirac structure $\lL$ is called
\emph{admissible} if there exists a smooth vector field $X_f$ on
$M$ such that $(X_f,df) \in \Gamma(\lL)$. We will denote by
$C^\infty_{\lL}(M)$ the set of $\lL$-admissible functions on $M$.
In the Poisson and symplectic cases the set of admissible
functions is all of  $C^\infty (M)$, but in general it is not the
case. If a function $f$ is admissible, we will call a vector field
$X_f$ such that $(X_f,df)$ is a section of $\lL$ a
\emph{Hamiltonian vector field} associated to $f$. In \cite{C} it
is shown that, in spite of the fact that Hamiltonian vector fields
are not unique in general, the bracket
\begin{equation}\label{E:PoissonBracket}
    \{ f, g \} = X_f(g)
\end{equation}
defines a Poisson algebra structure on the space
$C^\infty_{\lL}(M)$ of $\lL$-admissible functions on $M$ (see also
\cite{BW}).
\\ \\
The Courant bracket (\ref{E:CB}) can be {\em twisted} by an extra
term given by a $3$-form on $M$. In \cite{SW} it is pointed out
that brackets of the form
\begin{equation}\label{E:TwistedCB}
   [ X \oplus \xi, Y \oplus \eta ]_H
   = [ X \oplus \xi, Y \oplus \eta ]_C - \iota_Y \iota_X H,
\end{equation}
where $X \oplus \xi, Y \oplus \eta \in \Gamma (\T M)$ and $H \in
\Omega^3_{cl}(M)$ is a closed $3$-form on $M$, called the
twisting, give rise to the same kind of structure as before. A
maximal isotropic sub-bundle $\lL$  of $\T M = TM \oplus T^* M$,
which is involutive under the twisted Courant bracket
(\ref{E:TwistedCB}) on $\T M$ is called a {\em $H$-twisted Dirac
structure}. We will denote such structures by $\lL_H$ ---although
the twisting actually appears on the bracket and not on the
sub-bundle $\lL$--- to distinguish the twisted and non-twisted
cases. Twisted Dirac structures appear naturally in Poisson
geometry when, for example, a reduction of a (twisted or
non-twisted) Dirac structure is performed \cite{BCG}. In quantum
field theory and superstring theory, the form $H$ has an
interpretation as the Neveu-Schwarz $3$-form
\cite{Gra}. \\
\\
In this paper we address the question of the definition of Poisson
algebras of smooth functions on $M$ associated to {\em twisted}
Dirac structures (the non-twisted case has been studied in
\cite{C}). We will define in section \ref{S:Admissible} the set of
$H$-\emph{admissible functions} associated to a twisted Dirac
structure $\lL_H$, and we will show in theorem \ref{T:Theorem}
that this set has the structure of a Poisson algebra with the
usual Poisson bracket. Our definition of admissible function in
the twisted case was inspired by the notion of Hamiltonian
symmetries given in \cite{U}, in the context of differential
graded Lie algebras associated to dg-manifolds, so we will follow
such a point of view. Many results proven in the case of
admissible functions are also true in the case of admissible
pairs, in the sense of definition \ref{D:Admissible} (see also
\cite{BHR}\cite{BS}\cite{Z}), but here we want to focus on the
case of functions on $M$. In section \ref{S:Admissible} we prove
our main results and, in order to illustrate their significance in
the case of twisted Dirac structures defined by graphs of
non-degenerate $2$-forms, as in (\ref{E:SymplGraph}), we give a
characterization of our notion of $H$-admissibility in this
situation.  We end this paper by illustrating the non-triviality
of the algebra of $H$-{admissible functions} associated to a
twisted Dirac structure with a well-known example from the
classical theory of dynamical systems, which arises naturally in
this context.

\vspace{0.5cm}

\section{The Poisson algebra associated to twisted Dirac
structures}\label{S:Admissible}

In this section we will describe the Poisson algebra associated to
a $H$-twisted Dirac structure $\lL_H$ on a closed smooth manifold
$M$, namely the algebra $C^\infty_{\lL_H}(M)$ of $H$-admissible
functions associated to the Dirac structure. Even though many
facts we want to show explicitly for Dirac structures in $\Gamma
(TM \oplus \Lambda^{n-1}T^*M)$, when $n=2$, are true for any $n\ge
0$ (i.e. for any higher analogue in the sense of
\cite{BS}\cite{Z}), we will focus on this case, where the Poisson
algebras of {\em functions} associated to twisted Dirac structures
appear. Our definition of admissible function in the twisted case
was inspired by the notion of Hamiltonian symmetries given in
\cite{U}, in the context of differential graded Lie algebras
associated to dg-manifolds, so we will begin by a short recall of
such a point of view in order to motivate our approach.

\subsection{Derived Brackets and Hamiltonian
Symmetries}\label{S:HamiltonianSymmetries} Let us consider, for $n
\ge 0$, the trivial $\R[n]$-bundle $P_{[n]}$ over the odd tangent
bundle $T[1]M$, with the derivation given by
\begin{equation}\label{E:QH}
    Q_H = d + H \partial_t,
\end{equation}
where $d$ denotes the de Rham differential and $H \in
\Omega^{n+1}(M)$. The derivation $Q_H$ defines a homological
vector field (i.e. it satisfies $[Q_H, Q_H]=0$) if and only if $H$
is closed, and isomorphism classes of such bundles over $T[1]M$
are in one to one correspondence with $H^n_{dR} (M, \R)$ (see
\cite{Roy}). Smooth functions on the dg-manifold ($P_{[n]}, Q_H$)
correspond, locally, to the algebra $C^\infty (P_{[n]})=
\Omega^\bullet (M) \otimes S[t]$, where $t$ denotes the coordinate
on the fiber $\R$, whose degree is defined to be $n$, so that
$Q_H$ is actually a derivation of degree $1$ (see
\cite{CS}\cite{Roy}\cite{S}\cite{U} and \cite{V1} for the
background, original references and notations related with the
point of view of graded manifolds, and \cite{Va} for the
relation between homological vector fields and Lie algebroids).
\\
\\
The twisted Courant bracket (\ref{E:TwistedCB}) is known to be the
derived bracket obtained from the complex of derivations
$\textsf{Der}^\bullet(P_{[2]}, Q_H)$ of the DGLA (differential
graded Lie algebra) associated to the dg-manifold $(P_{[2]}= T[1]M
\oplus \R[2], Q_H)$ (see \cite{Roy}\cite{S}).  Since the derived
complex of derivations associated to the dg-manifold ($P_{[n]},
Q_H$) is nothing but the extended de Rham complex \cite{U}
\begin{equation}\label{E:DDer}
\begin{array}{ccccccccccccccccccc}
\Omega^{0}(M) \!\!\! & \!\stackrel{d}{\rightarrow} \!\!\! & \!\!\!
\Omega^{1}(M) \!\!\! & \! \stackrel{d}{\rightarrow} \!\!\! &
\!\cdots \Omega^{n-2}(M)\!\!\! \!\!\!\! & \! &
\stackrel{d}{\rightarrow} & \!\!\! \mathfrak{X}(M) \oplus
\Omega^{n-1}(M) ,
\end{array}
\end{equation}
and the corresponding derived brackets are
\begin{eqnarray*}
  \lc \alpha \partial_t , \beta \partial_t \rc &=& [[d + H \partial_t,
\alpha
\partial_t], \beta \partial_t] = [(d  \alpha)\partial_t , \beta \partial_t ] = 0 \;\;\;\;\;\;\;\;\;\;\;\;
\forall  \alpha, \beta \in \Omega^{n-k}(M),\\
  \lc i_X + \alpha \partial_t , \beta \partial_t \rc &=& [ \mathcal{L}_X + (d \alpha +
i_X H)\partial_t , \beta \partial_t ] =(\mathcal{L}_X \beta +
  i_X H)\partial_t
\;\;\;\;\;\;\;\;\;\;\; \forall \beta \in \Omega^{n-k}(M),
\end{eqnarray*}
and
\begin{equation}\label{E:DerBraTwist}
\lc i_X + \alpha \partial_t , i_Y +  \beta \partial_t \rc =
i_{[X,Y]} + ( \mathcal{L}_X \beta - i_Y d \alpha - i_Y i_X
H)\partial_t \;\;\;\;\;\;\;\;\;\;\;\;\;\;\;\;\; \forall \alpha,
\beta \in \Omega^{n-1}(M),
\end{equation}
it follows that, for any $n \ge 0$, we have a twisted Courant
bracket on sections of the bundle $TM \oplus \Lambda^{n-1} T^*M$
defined as the anti-symmetrization of (\ref{E:DerBraTwist}). A
\emph{Dirac structure of type $n \ge 1$}, is an isotropic
sub-bundle $\lL$ of $\Gamma(TM \oplus \Lambda^{n-1} T^*M) \cong
\frak{X}(M)\oplus \Omega^{n-1} (M)$ such that $\lL^\perp = \lL$,
with respect to the symmetric pairing
\begin{equation}\label{E:SPairingn}
    \langle X \oplus \alpha, Y \oplus \beta \rangle_{+} = {1 \over 2} (i_X \beta +i_Y \alpha),
\end{equation}
for $X \oplus \alpha, Y \oplus \beta \in \frak{X} (M) \oplus
\Omega^{n-1}(M)$, and such that $\lc \Gamma(\lL), \Gamma(\lL) \rc
\subset \Gamma(\lL)$ with respect to the  bracket defined by
(\ref{E:DerBraTwist}). These higher analogues of Dirac structures
and Courant algebroids have been studied in \cite{Z} and
\cite{BS}, respectively, and they are related with the $n$-plectic
structures defined in \cite{BHR}. It has been shown that, in the
non-twisted case (i.e. when the $(n+1)$-form $H$ is zero) there
are Poisson algebras of forms associated to them. When $n=2$,
these algebras correspond to the Poisson algebras of
\emph{admissible functions} on $M$, associated to Dirac structures
whenever $H=0$, defined by Courant in \cite{C}. If the twisting
$H\in \Omega^3 (M)$ is non-trivial, there is a stronger idea of
admissibility from which a Poisson structure can be given to a
subspace of $C^\infty (M)$. This idea of admissibility is encoded
in the sub-DGLA of derivations of $(P_{[n]}, Q_H)$ consisting of
infinitesimal symmetries of the bundle $P_{[n]}$ obtained from
``geometric symmetries" of the twisting form $H$, defined in
\cite{U}. In particular, the derivations in degree $-1$ in such a
complex are the given by,
\begin{equation}\label{E:Der0}
    {\textsf{GDer}}^{-1}(Q_H) = \left\{ i_X  + \alpha \partial_t  \in {Der}^{-1}(Q_H) \mid d\alpha+ i_X H=0
    \right\}
\end{equation}
and, since the condition $d\alpha+ i_X H=0$, for $\alpha \in
\Omega^{n-1}(M)$, ensures that $ \mathcal{L}_XH  =0$, it contains
the symmetries of $H$ encoded by vector fields on $M$. The derived
algebra of this sub-DGLA, called the \emph{Hamiltonian symmetries}
of the homological vector field $Q_H$ in \cite{U}, turns to be the
extended de Rham complex (\ref{E:DDer}) with derived brackets
\cite{U}
\begin{eqnarray*}
  \lc \alpha \partial_t , \beta \partial_t \rc &=&  0 \;\;\;\;\;\;\;\;\;\;\;\;\;\;\;\;\;\;
\forall  \alpha, \beta \in \Omega^{n-k}(M),\\
  \lc i_X + \alpha \partial_t , \beta \partial_t \rc &=& (\mathcal{L}_X \beta)\partial_t
  \;\;\;\;\;\;\;\; \forall \beta \in
\Omega^{n-k}(M),
\end{eqnarray*}
and
\begin{equation}\label{E:GDerBraTwist}
\lc i_X + \alpha \partial_t , i_Y +  \beta \partial_t \rc =
i_{[X,Y]} + ( \mathcal{L}_X \beta )\partial_t
\;\;\;\;\;\;\;\;\;\;\;\;\;\;\forall \alpha, \beta \in
\Omega^{n-1}(M).
\end{equation}
We will use this dg-Leibniz algebra to give a meaning to
admissible function, giving rise to a Poisson algebra of smooth
functions on $M$ in the twisted case.

\vspace{0.5cm}
\subsection{Admissible pairs associated to Dirac structures}

\begin{defin}\label{D:Admissible} Let us consider $X \in \frak{X}(M)$ and $\alpha \in
\Omega^{n-1}(M)$ such that $i_X + \alpha\partial_t \in
Der^{-1}(Q_H)$. We say that $(X, \alpha)$ is an admissible pair,
and then $X$ is called a Hamiltonian vector field associated to
$\alpha$, if $i_X + \alpha
\partial_t \in {\textsf{GDer}}^{-1}(Q_H)$, i.e. if
\begin{equation}\label{E:Admissible}
d\alpha+ i_X H=0,
\end{equation}
for the twisting form $H$.
\end{defin}
\begin{exam}\label{Example0}
The first non-trivial case, our motivating example, is given by a
twisting by a closed $2$-form $h$ (i.e. when $n=1$). In this case
the pairing (\ref{E:SPairingn}) is identically zero on any pair of
sections of $TM \oplus 1$, and  (\ref{E:Admissible}) says that a
function $f \in C^\infty (M)$ is admissible, and has Hamiltonian
vector field $X$, if
\begin{equation}\label{E:Admn1}
df + i_X h = 0.
\end{equation}
Thus, if $h$ is also non-degenerate (i.e. a symplectic form), for
every smooth function on $M$ there exists a vector field $X$ such
that $(X, f)$ is an admissible pair, namely the Hamiltonian vector
field given by (\ref{E:Admn1}), in agreement with the classical
setting of symplectic geometry. Notice that, restricting the
derived bracket (\ref{E:GDerBraTwist}) to pairs of Hamiltonian
vector fields and admissible functions, we find
$$\lc i_{X_f} + f \partial_t , i_{X_g} +  g \partial_t \rc =
i_{[X_f,X_g]} + ( \mathcal{L}_{X_f} g )\partial_t$$ and thus, if
$h$ is non-degenerate,
$$\lc i_{X_f} + f \partial_t , i_{X_g} +  g \partial_t \rc =
i_{[X_f,X_g]} +  \{ f,g\}  \partial_t,$$ where
$$  \{ f,g\} = \mathcal{L}_{X_f} (g) = X_f (g) $$
is the usual Poisson bracket on functions associated to the
symplectic form $h$.
\end{exam}
\begin{exam} Replacing the closed non-degenerate $2$-form
of the preceding example by a closed non-degenerate $(n+1)$-form
on $M$,  we naturally get the corresponding notion of Hamiltonian
forms, Hamiltonian vector fields and the so-called
\emph{Hemi-bracket} in the context of $n$-plectic or
multi-symplectic geometry studied in  \cite{BHR}.
\end{exam}
Let us compare our approach to admissible functions to the one
used for Dirac structures in the literature (see \cite{C}).
Consider the image under $1 \oplus d$ of the space of admissible
sections sections of $TM \oplus 1$ in $\frak{X} (M) \oplus
\Omega^1(M)$, i.e. the set of pairs $(X, df)$ such that
(\ref{E:Admn1}) follows for the closed form $h \in \Omega^2(M)$.
If $h$ is non-degenerate, such an image defines a
\emph{non-twisted} Dirac structure on $M$, which is actually the
graph (\ref{E:SymplGraph}) of the isomorphism induced by $h$
between tangent and cotangent spaces of $M$ point by point. Thus,
the image under the exterior differential of the symplectic model
at the level $n=1$ is the Dirac symplectic model
(\ref{E:SymplGraph}) at the level $n=2$ \emph{without twisting}.
In this case it is easy to see that, if we consider $\alpha = df$,
the condition (\ref{E:Admissible}) is empty, so that any exact
$1$-form is admissible and the bracket (\ref{E:PoissonBracket})
defines a Poisson bracket on $C^\infty (M)$, as proven in
\cite{C}. Moreover, restricting the derived bracket
(\ref{E:DerBraTwist}) to admissible pairs gives
\begin{equation}\label{E:NonTwistPoissRel}
\lc i_{X_f} + df \partial_t , i_{X_g} +  dg \partial_t \rc =
i_{[X_f,X_g]} +  d \{ f,g\}  \partial_t.
\end{equation}
We will see later that the same can be done in the twisted case,
asking not only $f$, but also $df$, to be admissible in the sense
of Definition \ref{D:Admissible}. These facts are also true for
any $n \ge 1$ (see also \cite{BS}\cite{Z}).
\begin{prop}
The image by the de Rham differential of the set of admissible
pairs $(X, \alpha)$ of a $H$-twisted Dirac structure in
$\frak{X}(M) \oplus \Omega^{n-1}(M)$, for $H$ non-degenerate,
defines a non-twisted Dirac structure in $\frak{X}(M) \oplus
\Omega^{n}(M)$.
\end{prop}

\emph{Proof.} Let $(X,\alpha), (Y,\beta)$ be an admissible pair in
$\frak{X} (M) \oplus \Omega^{n-1}(M)$, for a $H$-twisted Dirac
structure $\lL_H \le \Gamma (TM \oplus \Lambda^{n-1}T^*M)$. The
sections $(X,\alpha), (Y,\beta)$ belong to the Dirac structure
$\lL_H$, so that the pairing (\ref{E:SPairingn}) on $(X, \alpha),
(Y, \beta) \in \Gamma (TM \oplus \Lambda^{n-1}T^*M)$ is zero, and
$i_X \beta = - i_Y \alpha$. Since both pairs are admissible, we
have that $d \alpha = -i_X H$ and $d \beta = -i_Y H$. Then, the
sections $(X,d\alpha), (Y, d\beta)$ in the image under $1 \oplus
d$ of $\lL_H$ satisfy
$$ \langle (X,d\alpha), (Y,d\beta) \rangle_+ = i_X d \beta + i_Y d \alpha = 0.$$
On the other hand, computing the non-twisted Courant bracket,
\begin{eqnarray*}
 [(X,d\alpha), (Y,d\beta)]_C  &=&  ([X,Y], \mathcal{L}_X d\beta ) \\
    &=& ([X,Y], -\mathcal{L}_X i_Y H ) \\
    &=& ([X,Y], - i_{[X,Y]} H + i_Y  \mathcal{L}_X H )\\
    &=& ([X,Y], - i_{[X,Y]} H )
\end{eqnarray*}
so that $ [(X,-i_X H), (Y,-i_Y H)]_C  = ([X,Y], - i_{[X,Y]} H )$
and, provided $H$ is non-degenerate, the result is proven $\Box$\\
\\
Notice that, in general, for $(X,\alpha), (Y,\beta)$ admissible
pairs, we have a natural candidate to define the Poisson bracket
between $\alpha , \beta \in \Omega^{n-1}(M)$, namely
$\mathcal{L}_X(\beta)$; it is elementary to show that, on
admissible elements, this ``bracket" is skew-symmetric and has
other nice properties. Actually, the twisting by $H$ here
corresponds to the \emph{$n$-plectic structures} discussed in
\cite{BHR}, when $H$ is a non-degenerate form. We will now
concentrate in the case $n=2$, for $\alpha$ and $\beta$ exact
forms, i.e. we go back to define a Poisson algebra associated to a
twisted Dirac structures in $\Gamma (\T M)$.

\subsection{Admissible functions in the twisted case}

Let us now turn to the case of a $H$-twisted Dirac structures
$\lL_H$ in $\Gamma (\T M)$, where $H \in \Omega^3(M)$ is closed
(see \cite{BW} and \cite{SW}). Let us denote by $T_{\lL_H}$ the
tensor defined on sections of $\T M$ by
\begin{equation}\label{E:TensorL}
T_{\lL_H} (A \otimes B \otimes C) = \langle [ A, B ]_H , C
\rangle_+,
\end{equation}
where the pairing (\ref{E:SPairing}) and the twisted Courant
bracket (\ref{E:TwistedCB}) are used. This tensor, sometimes
called the Courant tensor, was defined in the non-twisted case in
\cite{C} in order to show that the bracket defined by
(\ref{E:PoissonBracket}), in the case of a non-twisted Dirac
structure $\lL$, defines a Poisson algebra on the algebra
$C^\infty_{\lL}(M)$ of admissible functions on $M$. Indeed, it is
clear from the definition that for any Dirac structure $\lL$ we
have $T_{\lL} (A \otimes B \otimes C) = 0$ whenever $A,B,C \in
\Gamma (\lL)$ and since, given admissible pairs $(X_f, df), (X_g,
dg)$ and $(X_h, dh) \in \Gamma (\lL)$, it is shown in \cite{C}
that
$$T_{\lL} ((X_f, df) \otimes (X_g, dg) \otimes(X_h, dh)) =
 \{f,\{g,h \}\} +  \{g,\{h, f \}\}  +   \{h,\{f,g \}\} ,$$
then Jacobi identity follows.
\begin{prop}\label{P:TLH}
Let $T_{\lL_H}$ denote the Courant tensor associated to a twisted
Dirac structure $\lL_H$ given by (\ref{E:TwistedCB}). Then
 $$T_{\lL_H} ((X_1, \alpha_1) \otimes (X_2, \alpha_2) \otimes(X_3, \alpha_3)) =
 T_{\lL} ((X_1, \alpha_1) \otimes (X_2, \alpha_2) \otimes(X_3, \alpha_3))
 - i_{X_3}i_{X_2}i_{X_1} H ,$$
for any $(X_1, \alpha_1), (X_2, \alpha_2), (X_3, \alpha_3) \in
\Gamma (\lL_H)$, where $T_{\lL}$ denotes the tensor
(\ref{E:TensorL}) associated to the non-twisted Courant bracket.
\end{prop}

\emph{Proof. } Let us take sections $A_i =(X_i, \alpha_i)$, for
$i=1,2,3$, in $\Gamma(\lL_H)$. By definition $$T_{\lL_H} (A_1
\otimes A_2 \otimes A_3) = \langle [ (X_1, \alpha_1), (X_2,
\alpha_2) ]_H , (X_3, \alpha_3) \rangle_+,$$ so that, from
(\ref{E:TwistedCB}),
$$T_{\lL_H} (A_1 \otimes
A_2 \otimes A_3) = \langle ([ X_1,X_2] , \mathcal{L}_{X_1}\alpha_2
-\mathcal{L}_{X_2} \alpha_1 + d (\alpha_1(X_2)) - i_{X_2}i_{X_1}H
), (X_3, \alpha_3) \rangle_+,$$ and then
\begin{eqnarray*}
 T_{\lL_H} (A_1 \otimes
A_2 \otimes A_3) &=& i_{X_3}(\mathcal{L}_{X_1}\alpha_2
-\mathcal{L}_{X_2} \alpha_1 + d
(\alpha_1(X_2)) - i_{X_2}i_{X_1}H) + i_{[ X_1,X_2]}\alpha_3 \\
   &=&  i_{X_3}(\mathcal{L}_{X_1}\alpha_2 -\mathcal{L}_{X_2} \alpha_1 + d
(\alpha_1(X_2)) ) + i_{[ X_1,X_2]}\alpha_3- i_{X_3}i_{X_2}i_{X_1}H\\
   &=&  T_{\lL} ((X_1, \alpha_1) \otimes
(X_2, \alpha_2) \otimes(X_3, \alpha_3)) - i_{X_3}i_{X_2}i_{X_1}H.
\end{eqnarray*}
 $\Box$
\\ \\
It follows that the Jacobi identity, for brackets
(\ref{E:PoissonBracket}) of admissible functions, has an
obstruction in the twisted case given by
$$ T_{\lL_H} ((X_f, df) \otimes (X_g, dg) \otimes(X_h, dh)) =
 \{f,\{g,h \}\} +  \{g,\{h, f \}\}  +   \{h,\{f,g \}\} - H (X_f, X_g , X_h) .$$
We will see that this obstruction disappears if we restrict the
space of admissible functions, as defined in \cite{C}, to a
smaller set in which the admissibility in the sense of definition
\ref{D:Admissible} plays a central role.

\begin{defin}
Let $\lL_H$ be a Dirac structure on a manifold $M$, twisted by a
$3$-form $H$. A function $f$ is $H$-admissible if it is admissible
in the sense of Courant and $(X_f, df)$ is an admissible pair in
$\Gamma (\lL_H)$. We will denote by $C^\infty_{\lL_H}(M)$ the set
of $H$-admissible functions on $M$.
\end{defin}

Notice that if $H=0$, i.e. if there is no twisting, the definition
of admissible function coincides with the one of Courant. On the
other hand, if the twisting is non-trivial, the set of
$H$-admissible functions may be smaller than the space of
admissible functions in the usual sense but, as we will see, it is
still a Poisson algebra. First, we show that we recover the usual
bracket relation (\ref{E:NonTwistPoissRel}), but this time with
$H$-admissible functions and the twisted bracket:
\begin{prop}\label{P:PoissBrakAdm}  Restricting the twisted Courant bracket
(\ref{E:TwistedCB}) to admissible pairs  $(X_f, df)$ and $(X_g,
dg)$ gives
$$[(X_f, df), (X_g, dg) ]_H =
([X_f,X_g] ,  d \{ f,g\}) .$$
\end{prop}

\emph{Proof.} Since $(X_f, df), (X_g, dg) \in \Gamma(\lL_h)$ are
$H$-admissible,
\begin{eqnarray*}
  [ (X_f, df),(X_g, dg)]_H &=&  ([X_f,X_g] , \mathcal{L}_{X_f} (dg) - i_{X_g} (d^2f + i_{X_f} H))\\
    &=&  ([X_f,X_g] , d \mathcal{L}_{X_f} (g)) \\
    &=& ([X_f,X_g] ,d \{ f,g\})
\end{eqnarray*}
$\Box$
\\

\begin{thm}\label{T:Theorem}
Let $f,g$ be $H$-admissible functions on $M$ with respect to the
twisted Dirac structure $\lL_H$, where $H \in \Omega^3(M)$ is
closed. Then the product $fg$ and the bracket $\{f,g\}$ defined by
(\ref{E:PoissonBracket}) are $H$-admissible functions. Moreover,
such a bracket satisfies both Leibniz and Jacobi identities, and
 defines a Poisson algebra structure on the space
$C^\infty_{\lL_H} (M)$.
\end{thm}

\emph{Proof.} Let $f,g,h$ be $H$-admissible functions, and let us
denote by $(X_f, df) , (X_g, dg)$, $(X_h, dh)\in \Gamma (\lL_H)$
the corresponding admissible pairs in ${\textsf{GDer}}^{-1}(Q_H)$.
Let $X_{fg}= gX_f + f X_g $, then
$$g(X_f, df) + f (X_g, dg)=
(gX_f + f X_g, gdf + fdg) = (X_{fg}, d(fg)) \in \Gamma (\lL_H)$$
and
$$ i_{X_{fg}}H =  g \, i_{X_f}H + f \, i_{X_g}H =0,$$
so that $fg$ is also $H$-admissible. Notice that both antisymmetry
and Leibniz identity are independent of the twisting. Indeed,
since $  i_{X_g} df = -  i_{X_f} dg$ for admissible pairs $(X_f,
df)$ and $ (X_g, dg)$,
$$  \{ f, g\}  = X_f (g) = \mathcal{L}_{X_f} (g)=
i_{X_f} dg + d \, i_{X_f} g =-  i_{X_g} df = - \mathcal{L}_{X_g} (f) = -\{ g, f\} $$
and, second,
$$ \{ fg, h\} = X_{fg}(h) = gX_f (h)+ f X_g(h) = g \{ f, h\} + f\{ g, h\}.$$
Next, since $$i_{[X_f,X_g] } H = \mathcal{L}_{X_f } i_{X_g} H -
i_{X_g}  \mathcal{L} _{X_f } H = - i_{X_g} d i_{X_f } H + i_{X_g}
i_{X_f } dH =0,$$ proposition \ref{P:PoissBrakAdm} implies that
$\{ f, g\} $ is $H$-admissible and
$$X_{\{ f, g\} } = [X_f,X_g] . $$ Finally,
definition \ref{D:Admissible} implies that $i_{X_f}H = i_{X_g}H
=i_{X_h}H =0$, so that by Proposition \ref{P:TLH} there is no
obstruction to the Jacobi identity $\Box$

\begin{exam} Let $G$ be a Lie group with Lie algebra
$\mathfrak{g}$, and let $\langle , \rangle_{\mathfrak{g}}$ be a
 non-degenerate symmetric bilinear form on it. Then
$$\lL_G = \{ ((X_R - X_L ), {1 \over 2} ( X_R - X_L ) ) |  X \in \mathfrak{g}\} \le TG \oplus T^*G ,$$
where $X_R$ and $X_L$ denote the right-invariant and
left-invariant vector fields associated to $X \in \mathfrak{g}$,
respectively, defines a twisted Dirac structure on $G$ with
twisting $3$-form
\begin{equation*}
    H_G = {1 \over 2} \langle [X,Y], Z  \rangle_{\mathfrak{g}},
\end{equation*}
called the Cartan-Dirac structure  on $G$ (see \cite{BW}). Since,
for an orthonormal basis $\{X_1, X_2 , \dots , X_n\}$ for
$\mathfrak{g}$,
$$ i_{X_l}i_{X_m}i_{X_n} H_G = \mathrm{C}^{lm}_n,$$
the structure constants of the Lie algebra, it is clear that the
space of $H_G$-admissible functions is completely determined by
$\mathfrak{g}$. In the case of $G=SO(3)$, with the usual
bi-invariant metric on it and the corresponding orthonormal basis
for $\mathfrak{so}(3)$, it is easy to see that $
i_{X_1}i_{X_2}i_{X_3} H_G =1$, so in this case the algebra of
$H_{SO(3)}$ admissible functions associated to the Cartan-Dirac
structure is trivial.
\end{exam}

\subsection{Twisted symplectic graphs and constants of motion} We will finish this paper
considering the example given by a symplectic graph twisted by a
closed $3$-form $H$ (also called $H$-closed $2$-forms in
\cite{SW}). Consider the Dirac structure defined in
(\ref{E:SymplGraph}), i.e. the graph $$\lL_h = \{(X, i_X h) \mid X
\in \frak{X}(M) \},$$ in $\T M$, of a non-degenerate $2$-form $h$.
It follows from the definition of the twisted bracket
(\ref{E:TwistedCB}) that this Dirac structure is integrable if and
only if $dh-H=0$, so that $h$ cannot be closed and, as a
consequence, the definition (\ref{E:PoissonBracket}) gives a
Poisson bracket for which, as follows from the remarks after
Proposition \ref{P:TLH} (see also \cite{SW}),
\begin{equation}\label{E:NoJac}
  \{f,\{g,h \}\} +  \{g,\{h, f \}\}  +   \{h,\{f,g \}\} =
  H (X_f, X_g , X_h) .
\end{equation}
Consider now functions  $f, g, h \in C^\infty_{\lL_H} (M)$ which
are $H$-admissible, then (\ref{E:Admissible}) implies that
$$i_{X_f}H = i_{X_g}H =  i_{X_h}H = 0,$$
so that the Jacobi identity holds. In this case, being a graph of
a symplectic form, the twisted Dirac structures associates to any
function a Hamiltonian vector field $X_f$, but it is this vector
field which makes $f$ an $H$-admissible function through the
condition $i_{X_f}H=0$. We can characterize such pairs with the
following
\begin{prop}\label{P:SymplGraph} A pair $(X_f, df)$ in $\lL_h$ is $H$-admissible if and
only if $\mathcal{L}_{X_f} h = 0$.
\end{prop}

\emph{Proof.} Since $({X_f}, df) \in \Gamma(\lL_h)$, $df=-i_{X_f}
h$ so that $\mathcal{L}_{X_f} h = d i_{X_f} h + i_{X_f} dh =
i_{X_f} H$, and the result follows $\Box$
\\ \\
Finally, as a consequence of Proposition \ref{P:PoissBrakAdm}, we
find back the usual bracket relation (\ref{E:NonTwistPoissRel})
for the $H$-twisted bracket:
$$  [ (X, df), (Y, dg) ]_H  = ([X_f,X_g] ,d \{ f,g\}).$$ \\
To see that the Poisson algebra of $H$-admissible functions is not
trivial in general, let $(M, \omega)$ be  a symplectic manifold
and consider the $2$-form $h= \varphi \cdot \omega$, where
$\varphi \in C^\infty (M)$. Then the twisted Dirac structure
(\ref{E:SymplGraph}) is integrable with respect to the twisted
Courant bracket (\ref{E:TwistedCB}) if and only if $H= dh = d
\varphi \wedge \omega$. Notice that, on the one hand, if $\varphi$
is chosen in such a way that $h$ is non-degenerate, any smooth
function on $M$ is admissible in the sense of Courant, i.e. for
any $f\in C^\infty (M)$ there exists a vector field $X_f \in
\mathfrak{X}(M)$ such that $df = - i_{X_f} h$, $(X_f, - df) \in
\Gamma (\lL_H)$,  where our notation means involutivity with
respect to the twisted Courant bracket (\ref{E:TwistedCB}). On the
other hand, by Proposition \ref{P:SymplGraph}, a smooth function
$f$ on $M$ is $H$-admissible if and only if
$$ \mathcal{L}_{X_f} h = (\mathcal{L}_{X_f}\varphi) \omega + \varphi ( \mathcal{L}_{X_f} \omega)
= \{f, \varphi \} \omega = 0.$$ This means that, in the cases in
which $\varphi$ is the Hamiltonian function for a dynamical system
with phase space $(M, \omega)$, an observable $f \in C^\infty (M)$
is $H$-admissible if and only if it is a constant of motion.
\begin{exam} \textbf{Angular Momentum.} Consider $M= T_{o}^*\R^3$ with canonical Darboux coordinates
$(q_1,q_2,q_3,p_1,p_2,p_3)$ and symplectic form $\omega =
\sum_{i=1}^3 dp_i \wedge dq_i$. Let us take $\varphi
(q_1,q_2,q_3,p_1,p_2,p_3)= \sum_{i=1}^3 {p_i^2 \over 2} + V(r)$,
where $r=(q_1^2 + q_2^2 + q_3^2)^{1 \over 2}$ and $V(r)$ denotes a
radial potential. Then the functions
\begin{eqnarray*}
  L_1 &=& q_2 p_3 - q_3 p_2 \\
  L_2 &=& q_3 p_1 - q_1 p_3 \\
  L_3 &=& q_1 p_2 -  q_2 p_1
\end{eqnarray*}
are $H$-admissible for $H = d \varphi \wedge \omega$, and
$$ \{ L_1, L_2\} = L_3. $$
\end{exam}

\vspace{0.9cm}

{\bf Acknowledgements.}  The author is grateful to Henrique
Bursztyn, Michel Cahen, Simone Gutt, Yoshiaki Maeda, Juan Camilo
Orduz, Bernardo Uribe, and Alain Weinstein for many stimulating
discussions on the geometry of Poisson and Dirac manifolds. The
author also thanks the referee for pointing him reference
\cite{BS}, and for his very pertinent comments which helped to
improve the exposition of this article. This research has been
supported by the \emph{Vicerrector\'ia de Investigaciones} and the
\emph{Faculty of Sciences} of the Universidad de los Andes.

\end{document}